\documentclass[11pt]{article}
\usepackage[margin=1in]{geometry}
\usepackage{setspace}
\usepackage{parskip}
\usepackage{amsmath}
\usepackage{amsfonts}
\usepackage{amssymb}
\usepackage{amsthm}
\usepackage{mathtools}
\usepackage{mathrsfs}
\usepackage{bm}
\usepackage{graphicx}
\usepackage{float}
\usepackage{subcaption}
\usepackage{caption}
\usepackage{booktabs}
\usepackage{multirow}
\usepackage{array}
\usepackage{tikz}
\usetikzlibrary{
    arrows.meta,
    calc,
    decorations.pathmorphing,
    positioning,
    shapes.geometric
}

\usepackage{tikz-cd}
\usepackage[
    colorlinks=true,
    linkcolor=blue,
    citecolor=blue,
    urlcolor=blue
]{hyperref}
\usepackage[numbers]{natbib}
\usepackage{enumitem}
\usepackage{xcolor}

\newtheorem{theorem}{Theorem}[section]
\newtheorem{lemma}[theorem]{Lemma}
\newtheorem{proposition}[theorem]{Proposition}

\theoremstyle{definition}
\newtheorem{definition}[theorem]{Definition}
\newtheorem{example}[theorem]{Example}
\newtheorem{remark}[theorem]{Remark}
\newcommand{\C}{\mathbb{C}}
\newcommand{\Ce}{\widehat{\C}}
\newcommand{\R}{\mathbb{R}}

\newcommand{\Z}{\mathbb{Z}}
\newcommand{\N}{\mathbb{N}}

\title{A Finite Automaton Translating Geometric Addresses into Dynamical Itineraries on the Basilica}
\author{Casey Ward, Chenzixi Zhao}
\date{\today}

\begin{document}

\maketitle

\begin{abstract}
Fractals are generally known for their self-similarity and repeating patterns. One way to produce fractals is via the dynamics of functions from the extended complex plane to itself. These fractals are known as Julia sets. In this paper, we explore how methods from symbolic dynamics and computation theory can be used to describe where points on the basilica Julia set might be located depending on their dynamical orbit under the polynomial function $f(z)=z^2-1$.
\end{abstract}

\tableofcontents

\newpage
\setcounter{section}{-1}
\section{Introduction}
Iterating a single quadratic polynomial already produces some of the most interesting sets in mathematics. In this paper, we study the Basilica, which is the Julia set of  $f(z)=z^2-1$, through two symbolic addresses of its points. The two addresses are the dynamical address, which records the itinerary of a point under $f$, and the geometric address, which records its position in a nested subdivision of the set. We determine which transitions can occur in a dynamical address and find a finite automaton that transfers geometric addresses into dynamical addresses.

\begin{definition}
Let $f(z)=z^2+c$ with $c\in\C$. The \textbf{filled Julia set} of $f$ is
\[
K_f=\left\{z\in\C:\sup_{n\ge 0}|f^n(z)|<\infty\right\},
\]
and the \textbf{Julia set} of $f$ is its boundary $J_f=\partial K_f$. The \textbf{Basilica} is the Julia set $J=J_f$ of $f(z)=z^2-1$.
\end{definition}

The graph of the Basilica is shown below. We will study its geometric and dynamical properties in the next section.

\begin{figure}[H]
    \centering
    \includegraphics[width=0.55\textwidth]{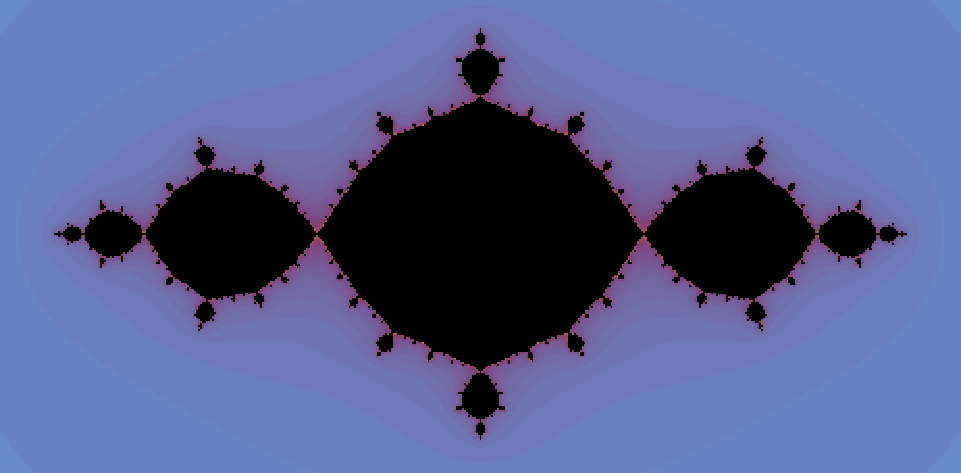}
    \caption{The Filled Basilica Julia Set}
\end{figure}



Intuitively, the Julia set describes the ``chaos'' of the orbit of $f$. To code the points of $J$ by symbols, we come up with two different ways: one based on the dynamics of the Basilica and another based on its
geometric appearance.

The \textbf{geometric address} comes from a different source. Following Belk and Forrest \cite{belk-forrest}, the Basilica is the limit space of an \emph{edge replacement system}: starting from a graph with four edges $t,l,b,r$, each edge is repeatedly replaced by three edges labeled $1,2,3$. A point of $J$ is then described by a sequence in $\{t,l,b,r\}\times\{1,2,3\}^{\N}$, read like the ternary expansion that locates a point of the Cantor set. Our second main result is using an \textbf{automaton machine} to translate between the two codings.

To do this, we construct a Mealy machine $M_D$ with five states that reads a geometric address one symbol at a time and writes the corresponding dynamical symbols. We show that $M_D$ is compatible with the dynamics of $f$, so the orbit of a point can be read off from its geometric address.

\section{Preliminaries in Julia Sets}

\subsection{Complex Dynamical Systems}
For the sake of discussion, we will restrict the family of maps in question to rational functions on the extended complex plane. With this in mind, we will begin by defining what a dynamical system is in the context of the Basilica.
\begin{definition}

For a complex dynamical system given by a map $f:\widehat{\C}\to\widehat{\C}$, the \textbf{Forward Orbit} of $z \in \mathbb{C}$ is:
\[
\mathcal{O}^{+}(z)=\{z,f(z),f^{2}(z),f^{3}(z),\ldots\},
\]
where $f^{n}(z)$ represents n-times iteration of $z$ under $f$ and $\Ce=\mathbb{C}\cup\{\infty\}$.
\end{definition}

In this paper, we will be focusing on the orbit of \emph{quadratic polynomials} with independent variables in the extended complex plane. 



\begin{definition} 
    Let $f:\Ce\to \Ce$ be a function. Then the forward orbit generated by a point $z_0$ under the dynamics of $f$ is said to be periodic if there is a natural number $n \in \N$ such that $f^{n}(z_0)=z_0$. The period of the orbit $t$ is defined as follows:
        $$t:=\min_{n\in \N}\{f^{n}(z_0)=z_0\}$$
\end{definition}

Another way that one could interpret $t$ is the smallest number for which $z_0$ is a fixed point of the function $f^t$.

\subsection{Böttcher Coordinate}
The exterior dynamics of the Basilica can be described using the
\textbf{B\"ottcher coordinate}. Since the critical orbit is bounded, $K$ is connected \cite{Milnor}, so there is a unique
conformal isomorphism
\[
\Phi\colon \C\setminus K\longrightarrow\{w\in\C:|w|>1\}
\]
normalized by $\Phi(z)/z\to 1$ as $z\to\infty$. It conjugates $f$ to
squaring, so $\Phi(f(z))=\Phi(z)^2$ for $z\notin K$ (see \cite{Milnor}).

For $\theta\in\R/\Z$, the \textbf{external ray} of angle $\theta$ is
\[
R_\theta=\{\Phi^{-1}(re^{2\pi i\theta}):r>1\}.
\]
Because the Basilica is locally connected \cite{Milnor}, every external ray lands on
$J$, and the landing map $\ell(\theta)=\lim_{r\downarrow 1}\Phi^{-1}(re^{2\pi i\theta})$
is continuous and onto. The conjugacy gives
\[
f(R_\theta)=R_{2\theta},\qquad f(\ell(\theta))=\ell(2\theta),
\]
where angles are taken modulo $1$. Thus iteration on the Basilica can be
studied through angle doubling on the circle.

The angles $1/3$ and $2/3$ form a two-cycle under doubling. By the
periodic ray landing theorem \cite{Milnor}, the rays landing at a
repelling fixed point all have the same period. The only points of $J$
fixed by $f^2$ are $-\varphi$ and $\beta$, since the other two solutions
of $f^2(z)=z$ are $0$ and $-1$. The fixed ray $R_0$ lands at $\beta$,
so $R_{1/3}$ and $R_{2/3}$ must land at $-\varphi$. These are the only
angles of period two, so they are the only rays landing there. Since
$f^{-1}(-\varphi)=\{-\varphi,\varphi\}$, taking the other two preimages
under angle doubling shows that $R_{1/6}$ and $R_{5/6}$ are exactly the
rays landing at $\varphi$.

These four angles divide $\R/\Z$ into the open arcs
\[
\begin{aligned}
I_T&=(1/6,1/3), & I_L&=(1/3,2/3),\\
I_B&=(2/3,5/6), & I_R&=(5/6,7/6).
\end{aligned}
\]
The four pieces of the Basilica are $J_X=\ell(I_X)$ for
$X\in\{T,L,B,R\}$. Indeed, the two rays at each cut point form a
Jordan curve when joined at their common landing point and at infinity.
These curves separate the left and right pieces, respectively, from the
rest. The rays at $1/6$ and $1/3$, joined at infinity and by the segment
$[-\varphi,\varphi]$, separate the top and bottom pieces. The interior
of that segment lies in the basin of $\{0,-1\}$, as
$|f^2(x)|<|x|$ for $0<|x|<\varphi$. Thus the separating curves meet $J$
only at the cut points. The landing map sends each angle sector to its
corresponding side away from those points. Each $\ell(I_X)$ is connected,
and the four images are exactly the components of
$J\setminus\{-\varphi,\varphi\}$.

For the transition diagram we can therefore compute
$f(J_X)=\ell(2I_X)$, where $2I_X$ denotes the image of $I_X$ under
angle doubling. The angles whose iterates hit a cut angle form a
countable set, since they are successive preimages of four angles.

\section{Symbolic Dynamics of the Basilica}
The fact that the basilica is a quadratic map makes it special as it only has at most three unique the Basilica has some basic properties points allows us to characterize the dynamics in a somewhat interesting way. In particular we will discover that the basilica may be partitioned into finitely many pieces. The fact that the Basilica can be divided into finitely many pieces might lend itself to a closer study of the set as a whole. Periodicity lends to promising prospects as the finite partition and boundedness of the Basilica lend a potential shift in perspective through symbolic dynamics. In this section we will explore two different approaches to tackling a symbolic representation of the Basilica.

\subsection{The Division of the Basilica}
Before we talk about the symbolic dynamics of the Basilica, we need to divide it into different parts based 
on geometry. Let's look at the following graph. We divide the Basilica into four pieces by four colors. Then we assign
the red, purple, blue, and green parts with $T,L,B,R$, which obviously means \textbf{Top, Left, Bottom} and \textbf{Right}. 
We call this the \textbf{Partition of the Basilica}.

\begin{figure}[H]
    \centering
    \includegraphics[width=0.55\textwidth]{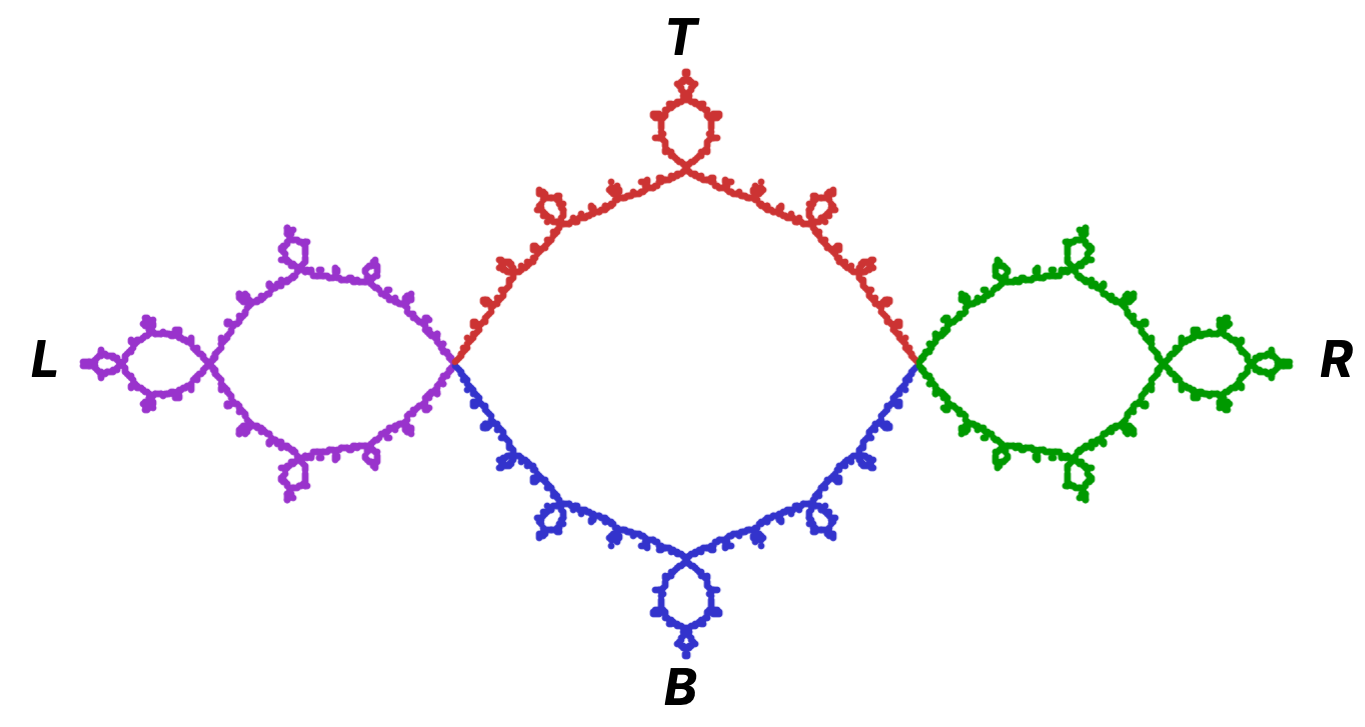}
    \caption{Partition of the Basilica}
    \label{fig:partition}
\end{figure}

There are some basic properties of the Basilica. It is important to know the actual coordinates of some points on the 
Basilica for future proofs.

\begin{proposition}[Tips, branch-cut points, symmetry]\label{prop:basilica}
\leavevmode
\begin{enumerate}
\item[(i)] \emph{(Branch-cut points.)} $-\varphi$ is a repelling fixed point of $f$, $f(\varphi)=-\varphi$, and $f^{-1}(-\varphi)=\{\varphi,-\varphi\}$. The points $\pm\varphi$ are the two branch-cut points at the waist of $J$.
\item[(ii)] \emph{(Tips and symmetry.)} $f(\beta)=\beta$ and $f(-\beta)=\beta$; the points $\pm\beta$ are the extreme points of $J$ on the real axis. Moreover, $J=-J=\overline{J}$.
\item[(iii)] \emph{(Partition by cuts.)} Consider the two vertical lines
\[
\ell_-=\{z\in\mathbb C:\operatorname{Re}z=-\varphi\},\qquad
\ell_+=\{z\in\mathbb C:\operatorname{Re}z=\varphi\},
\]
and the horizontal segment
\[
I=\{z\in\mathbb C:-\varphi<\operatorname{Re}z<\varphi,\ \operatorname{Im}z=0\}.
\]
Their complement has four open regions:
\[
\begin{aligned}
U_T&=\{z:-\varphi<\operatorname{Re}z<\varphi,\ \operatorname{Im}z>0\},\\
U_L&=\{z:\operatorname{Re}z<-\varphi\},\\
U_B&=\{z:-\varphi<\operatorname{Re}z<\varphi,\ \operatorname{Im}z<0\},\\
U_R&=\{z:\operatorname{Re}z>\varphi\}.
\end{aligned}
\]
Define $J_X=J\cap U_X$ for $X\in\{T,L,B,R\}$.
Using the geometric intersection property
\begin{equation}\label{eq:basilica-cut-intersections}
J\cap\ell_-=\{-\varphi\},\qquad
J\cap\ell_+=\{\varphi\},\qquad J\cap I=\varnothing,
\end{equation}
we obtain
\[
J\setminus\{-\varphi,\varphi\}=J_T\sqcup J_L\sqcup J_B\sqcup J_R.
\]
These are the four connected components of the punctured Basilica, and each $J_X$ is open in the relative topology of $J$.
\end{enumerate}
\end{proposition}

\begin{proof}
Since $\varphi^2=1-\varphi$ and $\beta^2=\beta+1$, direct substitution gives
\[
f(-\varphi)=-\varphi,\quad f(\varphi)=-\varphi,\quad f(\beta)=\beta,\quad f(-\beta)=\beta.
\]
Also $|f'(-\varphi)|=2\varphi>1$, so $-\varphi$ is repelling. Solving $z^2-1=-\varphi$ gives $z^2=\varphi^2$, hence the stated preimage pair. The repelling fixed points $-\varphi$ and $\beta$ belong to $J$; complete invariance then gives $\varphi,-\beta\in J$.
If $|z|>\beta$, then $|f(z)|\geq |z|^2-1>|z|$. The successive moduli increase and tend to infinity: a finite limit $L>\beta$ would satisfy $L\geq L^2-1$, which is impossible. Thus $K\subseteq\{|z|\leq\beta\}$, proving the assertion about the real extremes.
Since $f(-z)=f(z)$ and $f(\overline z)=\overline{f(z)}$, boundedness of the orbit is preserved by negation and conjugation. These homeomorphisms therefore preserve both $K$ and its boundary $J$.

For the geometric part, we first verify \eqref{eq:basilica-cut-intersections}. For real $x$ with $|x|\le r<\varphi$ we have $f^2(x)=-x^2(2-x^2)$, hence
\[
|f^2(x)|=|x|\cdot|x|(2-x^2)\le c\,|x|,\qquad c=r(2-r^2).
\]
Since $t\mapsto t(2-t^2)$ is increasing on $[0,1]$ and equals $\varphi(1+\varphi)=1$ at $t=\varphi$, we get $c<1$. Thus $f^{2n}(x)\to 0$, so $x$ lies in the basin of the superattracting cycle $\{0,-1\}$, which is open and contained in the Fatou set. Hence $J\cap I=\varnothing$. The equality $J\cap\ell_-=\{-\varphi\}$ follows from the identity $J\cap\{\operatorname{Re}z\le-\varphi\}=J_{(3)}$ in \cite[Section~6]{rogers-teplyaev}, since the cells there meet only at their boundary points; then $J\cap\ell_+=\{\varphi\}$ by $J=-J$.

By \eqref{eq:basilica-cut-intersections}, $J\setminus(\ell_-\cup\ell_+\cup I)=J\setminus\{\pm\varphi\}$, so $J\setminus\{\pm\varphi\}=J_T\sqcup J_L\sqcup J_B\sqcup J_R$. Each $U_X$ is open in $\mathbb C$, so $J_X=J\cap U_X$ is open in $J$ by the definition of the relative topology.

Each $J_X$ is nonempty. Indeed, $\beta\in J_R$ and $-\beta\in J_L$. Since $0\in\operatorname{int}K$ and $K$ is bounded, the ray $\{it:t>0\}$ leaves $K$ and therefore meets $J=\partial K$ at some $it_0$ with $t_0>0$; then $it_0\in J_T$, and $-it_0=\overline{it_0}\in J_B$ because $J=\overline{J}$.

Finally, every connected component of $J\setminus\{\pm\varphi\}$ lies in exactly one $J_X$, since the $J_X$ are disjoint, relatively open and cover $J\setminus\{\pm\varphi\}$. Hence each $J_X$ is a nonempty union of components. By \cite[Section~2]{rogers-teplyaev}, $J\setminus\{\pm\varphi\}$ has exactly four components, so each $J_X$ is exactly one component.
These components are the pieces $\ell(I_X)$ of Section~1.2, with the same labels. Indeed, $\beta$ and $-\beta$ are the landing points of the angles $0\in I_R$ and $1/2\in I_L$. Moreover, $\R\setminus K=R_0\cup R_{1/2}$, so the rays with angles in $(0,1/2)$ cannot cross the real axis and stay in the upper half-plane; hence $\ell(I_T)$ is the upper piece $J_T$, and $\ell(I_B)=J_B$ by conjugation.
\end{proof}

\begin{remark}
The identity $J\cap\{\operatorname{Re}z\le-\varphi\}=J_{(3)}$, and hence $J\cap\ell_\pm=\{\pm\varphi\}$, is used in \cite[Section~6]{rogers-teplyaev} without a separate proof; a numerical distance-estimate check along $\operatorname{Re}z=\pm\varphi$ is consistent with it.
\end{remark}
\subsection{Dynamic Address}

\begin{definition}
    Let $z \in J$ be a point on the Basilica. The \textbf{dynamical address} of $z$ is the string $(d_0 d_1 d_2 \cdots)$, 
    where $d_n \in \{T, L, B, R\}$ indicates which region of the partition contains $f^n(z)$.
\end{definition}

We will give some examples after we discuss the transitions of different addresses on the Basilica. 

\begin{theorem}\label{thm:transition}
    The following transition diagram 
lists all admissible 
symbolic transitions under the iteration:

\begin{figure}[H]
    \centering
    \includegraphics[width=0.55\textwidth]{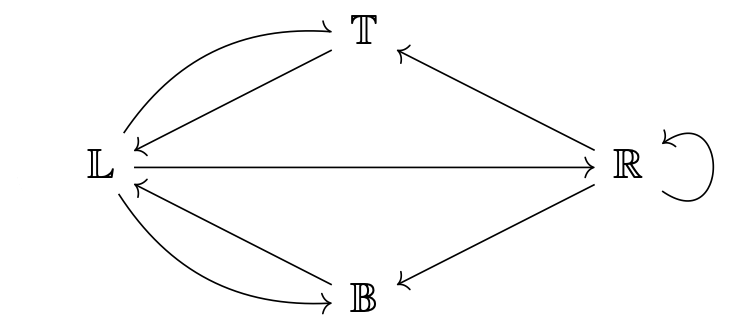}
    \caption{The Transition Diagram of the Basilica}
    \label{fig:transition}
\end{figure}
\end{theorem}

This transition diagram shows that for any point $z_0 \in \{T, L, B, R\}$, the new point $f(z_0)=z_0^2-1$ after one time of iteration
can only be located in the area where the arrows are pointing. We will start the proof with the following proposition.

We record the geometric facts about $J$ that will be used freely. The facts are stated in Proposition~\ref{prop:basilica}: they are the tips and the branch-cut
points visible in Figure~\ref{fig:partition}.

\subsection{Proof of Theorem~\ref{thm:transition}}

\begin{proof}
Let $\tau(\theta)=2\theta\pmod 1$. The landing relation gives
$f(J_X)=\ell(\tau(I_X))$. Directly doubling the four arcs gives
\[
\tau(I_T)=\tau(I_B)=I_L,\qquad
\tau(I_L)=\tau(I_R)=(2/3,4/3).
\]
The arc $(2/3,4/3)$ consists of $I_B$, $I_R$, $I_T$, and the two
angles $5/6$ and $1/6$ that land at $\varphi$. Hence
\[
f(J_T)=f(J_B)=J_L,\qquad
f(J_L)=f(J_R)=J_T\cup J_B\cup J_R\cup\{\varphi\}.
\]
This gives exactly the arrows in the diagram. Every listed arrow is
realized by an angle in a nonempty open subarc whose forward orbit
avoids the cut angles, since only countably many angles eventually
hit one of them.
\end{proof}

\subsubsection{Symbol Space}
Knowing that we have a partition of the basilica into finitely many pieces allows us to consider a new topological space $\Gamma_{D}$ which gives a symbolic description of the dynamics.

$$\Gamma_{D}=\{T,L,B,R\}^{\N},$$
with the product topology, in which two sequences are close when they agree on a long initial block. Not every sequence in $\Gamma_D$ is the address of a point: we call $\omega\in\Gamma_D$ \textbf{admissible} if every transition $\omega(n)\to\omega(n+1)$ appears in Theorem~\ref{thm:transition}, and we let $\Sigma_D\subset\Gamma_D$ be the closed set of admissible sequences. With the partition of the basilica we have, we may construct a quotient map $\pi_{D}:\Sigma_{D}\to J_{f}$ as follows. Write $\tau(\theta)=2\theta \pmod 1$ and $\overline{I_X}$ for the closed arcs of Section~1.2. For an admissible $\omega$ and $n\ge0$, let
$$I_{\omega(0)\cdots\omega(n)}=\{\theta:\ \tau^k(\theta)\in\overline{I_{\omega(k)}}\ \text{for}\ 0\le k\le n\}.$$
Since $\tau$ is injective on each $\overline{I_X}$ and maps it onto the union of the closed arcs of its successors, induction shows that $I_{\omega(0)\cdots\omega(n)}$ is a nonempty closed arc of length $2^{-n}|I_{\omega(n)}|\le 2^{-n}/3$. Hence these nested arcs shrink to a single angle $\theta(\omega)$, and we define
$$\pi_{D}(\omega)=\ell(\theta(\omega)).$$
Then $f^n(\pi_D(\omega))\in\overline{J_{\omega(n)}}$ for every $n\ge0$. If the orbit of $z\in J$ avoids $\pm\varphi$, then $z$ is the landing point of exactly one angle, its dynamical address is admissible, and $\pi_D$ sends this address to $z$. Two admissible sequences that agree up to $n$ give angles at distance at most $2^{-n}/3$, so $\theta$ is continuous, and so is $\pi_D=\ell\circ\theta$. Every angle has an admissible itinerary (at the countably many angles that hit a cut angle, take one-sided limits), so $\pi_D$ is onto. Because $\Sigma_{D}$ is compact and $J_f$ is Hausdorff, $\pi_D$ is a closed map, and so it is a quotient map. \qed

\subsection{Geometric Address}
\subsubsection{Geometric Address of the Basilica}
Different from the dynamical address, the geometric address is based on the geometric location of each 
point on the Basilica. The idea is similar to the itinerary of the Cantor set. To define the geometric 
address, we need to introduce two things first: the initial setting and the replacement rule.

We may still use the partition of the Basilica for the initial address, but in order to distinguish for dynamical address, 
we will change the four letters into lower case: 
$\{t,l,b,r\}$.


The idea of the replacement rule comes from James Belk and Bradley Forrest \cite{belk-forrest} and we can 
check the following graph for visualization. Without loss of generality, we choose a single edge to 
illustrate the rule first. Every time we apply the replacement, we take the midpoint of the edge and 
promote it to a new vertex; from this vertex we grow two edges back toward the endpoints of the original 
edge together with a loop, giving three new edges. We label these three edges $1,2,3$ in counterclockwise 
order, which gives a consistent orientation that is preserved under further iteration. Following the 
framework of Belk and Forrest \cite{belk-forrest}, we phrase this as an edge replacement system.

\begin{definition}
The \textbf{Basilica replacement system} is the edge replacement system 
$\mathcal{R} = (G_0,\, e \to R)$, where:
\begin{itemize}
    \item the \emph{base graph} $G_0$ consists of two vertices, two directed edges $t,b$ joining them, and 
    one loop at each vertex, $l$ and $r$; we identify its edge set with the four regions of the partition,
    \[
    E(G_0) = \{t,\, l,\, b,\, r\};
    \]
    \item the \emph{replacement rule} $e \to R$ replaces a single edge by promoting its midpoint to a new 
    vertex, from which three new edges emanate: two edges back to the endpoints of the original edge 
    together with a loop. These are labeled in counterclockwise order,
    \[
    E(R) = \{1,\, 2,\, 3\}.
    \]
\end{itemize}
As in \cite{belk-forrest}, the \emph{full expansion sequence} $\{G_n\}_{n=0}^{\infty}$ is obtained by applying 
$e \to R$ simultaneously to every edge of $G_n$ to form $G_{n+1}$, and the edges of $G_n$ are precisely 
the elements of
\[
E(G_0) \times E(R)^n = \{t,l,b,r\} \times \{1,2,3\}^n.
\]
\end{definition}

\begin{lemma}
The limit space of the Basilica replacement system $\mathcal{R}$ is homeomorphic to the Basilica $J$.
\end{lemma}

\begin{proof}
The base graph $G_0$ and replacement rule $e \to R$ described above form (an isomorphic copy of) the 
Basilica replacement system studied by Belk and Forrest. By their construction of the limit space of an 
edge replacement system \cite{belk-forrest}, the limit space of $\mathcal{R}$ is homeomorphic to the 
Julia set of $f(z) = z^2 - 1$. Since $J$ is exactly this Julia set, the limit space of $\mathcal{R}$ is 
homeomorphic to $J$.
\end{proof}

Now let us check some examples of the geometric address. 
\begin{definition}\label{def:geometric-address}
For a closed arc $[x,x+\lambda]$ of $\R/\Z$, lifted to $\R$, define its three
\emph{children}, labeled in increasing angular order, by
\[
[x,x+\lambda]_1=[x,\,x+\tfrac{\lambda}{4}],\qquad
[x,x+\lambda]_2=[x+\tfrac{\lambda}{4},\,x+\tfrac{3\lambda}{4}],\qquad
[x,x+\lambda]_3=[x+\tfrac{3\lambda}{4},\,x+\lambda].
\]
Starting from the closed initial arcs
\[
I_t=[1/6,1/3],\qquad I_l=[1/3,2/3],\qquad I_b=[2/3,5/6],\qquad I_r=[5/6,7/6]\pmod 1,
\]
that is, the closures of $I_T,I_L,I_B,I_R$, this defines recursively the arcs
$I_{s i_1\cdots i_n}=(I_{s i_1\cdots i_{n-1}})_{i_n}$ for $s\in\{t,l,b,r\}$ and
$i_k\in\{1,2,3\}$. Their lengths are
\[
|I_{s i_1\cdots i_n}| = |I_s|\prod_{k=1}^{n}\rho_{i_k},\qquad \rho=(\rho_1,\rho_2,\rho_3)=\bigl(\tfrac14,\tfrac12,\tfrac14\bigr),
\]
so the nested arcs shrink to a single angle. The corresponding \textbf{geometric cells} are
$C_{s i_1\cdots i_n}=\ell(I_{s i_1\cdots i_n})\subset J$, where $\ell$ is the landing map.
A \textbf{geometric address} of $z\in J$ is a sequence
\[
s\,i_1 i_2 i_3\cdots\in\{t,l,b,r\}\times\{1,2,3\}^{\N}
\qquad\text{with}\qquad z\in C_{s i_1\cdots i_n}\ \text{for every } n\ge 0,
\]
equivalently $z=\ell(\theta)$ where $\{\theta\}=\bigcap_{n} I_{s i_1\cdots i_n}$.
\end{definition}
\section{The Transfer of Address}
Up to this point, we have introduced two distinct symbolic descriptions associated with the Basilica. It is natural to ask whether there is a way to transfer the Dynamic Address and the Geometric
Address. The automaton machine is a good way we thought of that can deal with the issue.


To make this correspondence explicit, it is useful to view the translation
process as a symbolic input--output procedure, which is a sequence describing a
geometric address is processed according to fixed local rules to produce the
corresponding dynamical address. From this perspective, the translation can be
modeled by a finite automaton, which reads symbols sequentially and outputs the
appropriate dynamical symbols dictated by the underlying dynamics. An example
of this input--output structure is illustrated in the figure below:

\begin{figure}[H]
\centering

\begin{tikzpicture}[
  font=\small,
  x=0.9cm, y=1cm
]
\def\n{7} 
\def\inputbits{1,0,0,1,0,1,1}
\def\outputbits{0,1,0,0,1,1,0}

\def\halfheight{0.28} 
\node[anchor=east] at (-0.4, 0) {Input};
\node[anchor=east] at (-0.4,-1) {Output};
\draw (0,-\halfheight) rectangle (\n,\halfheight);
\foreach \k in {1,...,6} {
  \draw (\k,-\halfheight) -- (\k,\halfheight);
}
\foreach \b [count=\idx from 1] in \inputbits {
  \node at ({\idx-0.5},0) {\b};
}
\node[anchor=west] at (\n+0.3,0) {$\cdots$};

\begin{scope}[shift={(0,-1)}]
  \draw (0,-\halfheight) rectangle (\n,\halfheight);
  \foreach \k in {1,...,6} {
    \draw (\k,-\halfheight) -- (\k,\halfheight);
  }
  \foreach \b [count=\idx from 1] in \outputbits {
    \node at ({\idx-0.5},0) {\b};
  }
  \node[anchor=west] at (\n+0.3,0) {$\cdots$};
\end{scope}

\end{tikzpicture}

\caption{Example of input and output of an automaton}
\end{figure}

To achieve this, we require a formal structure capable of processing symbolic data in a controlled and deterministic way. 
This leads naturally to tools from automata theory, which provide a convenient language for encoding finite combinatorial rules 
governing symbolic transitions.

\subsection{The Automaton Machine}
We recall that an automaton consists of finitely many states together with
transition rules driven by symbols from a fixed alphabet. Since we must not
merely recognize valid sequences but convert an address in one system into an
address in another, we require a machine that reads and writes simultaneously:
a \emph{Mealy machine}. Its finiteness reflects the fact that the local
combinatorial structure of the Julia set is governed by finitely many
configurations, while its deterministic transition and output rules make the
translation well-defined across iterations.


  \begin{figure}[H]
  \centering
  \begin{minipage}{0.62\textwidth}
    \centering
    \includegraphics[width=\textwidth]{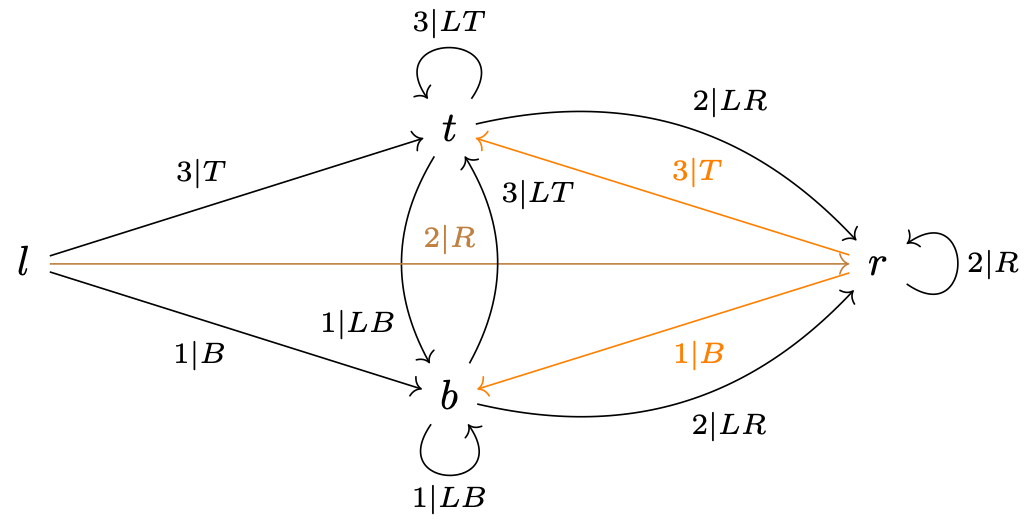}
    \caption*{Main automaton}
  \end{minipage}\hfill
  \begin{minipage}{0.32\textwidth}
    \centering
    \includegraphics[width=\textwidth]{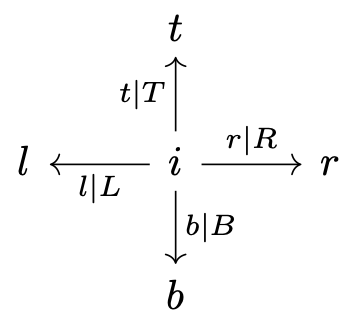}
    \caption*{Initial states}
  \end{minipage}
  \caption{The automaton machine translating geometric addresses into dynamical addresses}
  \label{fig:machine-gtd}
\end{figure}

The Mealy machine described above is an informal construction of what the potential machine could look like. Strictly speaking, the machine above is not a classical Mealy machine, as some of its arrows write two characters at once. In the construction of $M_D$ below we formalize it by allowing the output of each step to be a nonempty word.

The states $t,l,b,r$ record the geometric location of a point in the Basilica,
and the output alphabet $\{T,L,B,R\}$ encodes the dynamical address. Each
initial geometric symbol determines both a first output symbol and a starting
state, via the initial-state diagram on the right of
Figure~\ref{fig:machine-gtd}; the machine on the left then does the rest.

After this initialization, the machine reads the geometric address
sequentially. At each step, the current state reads one geometric symbol from
the input string. Among the outgoing arrows from that state, the arrow labeled
with the corresponding input symbol is selected. Each arrow is labeled in the
form
\[
\text{input} \mid \text{output},
\]
where the symbol to the left of the vertical bar denotes the geometric symbol
being read, and the symbol or symbols to the right denote the dynamical output
produced at that step. The arrow also determines the next state of the machine.
The following exploration will make clear why we would want to consider such a construction.
\begin{example}
 Consider the geometric address
\[
t323\ldots
\]
The initial symbol $t$ is interpreted using the initial-state diagram, which
outputs $T$ and places the machine in state $t$ in the main automaton. Reading
the next symbol $3$, one follows the arrow leaving state $t$ labeled $3\mid LT$,
producing the output $LT$ and returning to state $t$. The next symbol $2$ is
then read, and the arrow labeled $2\mid LR$ is followed, yielding the output
$LR$ and moving the machine to state $r$. Reading the next symbol $3$, the arrow
labeled $3\mid T$ is selected, producing the output $T$ and returning to the
appropriate state. Continuing in this manner, the geometric address is
converted into the dynamical address
\[
TLTLRT\ldots
\]

\end{example}

In the following section we formally define several variants of generalized sequence machines for finite state ``systems'' and construct a Mealy machine that relates the geometric address of a point to its dynamical orbit under $f$.

\subsection{Construction of $M_D$}

\begin{definition}\label{ref:DFA}
Let $\Lambda$ and $Q$ be finite sets.
A deterministic finite state automaton is a tuple $(Q,\Lambda,\alpha,q_0,F)$ consisting of a state set $Q$, an alphabet $\Lambda$, a transition function $\alpha:Q\times \Lambda\to Q$, an initial state $q_{0}\in Q$ and final states $F\subset Q$.
\end{definition}

The function $\alpha$ is defined inductively in the following way:

Let $\omega$ be a word in $\Lambda^{*}=\bigcup_{n\in\N}\Lambda^{n}$, with $\Lambda^{n}$ being the set of all words of length $n$ consisting of elements of $\Lambda$.

It may be the case that the automaton is only able to run when reading specific words in $\Lambda^{*}$. For this reason, we define the notion of acceptance.

\begin{definition}
    Let $M$ be a finite state automaton. We say a word $w\in \Lambda^{*}$ is accepted if $\alpha(q,w)\in F$. The set of all accepted words of an automaton is known as the language associated with $ M$, denoted $L(M)$.
\end{definition}

\begin{remark}
    Later we will use this fact to make an automaton whose language contains the set of all geometric addresses of $J_f$.
\end{remark}

$$\alpha(q_0,\omega_1)=q_j$$
$$\alpha(q_j,\omega_{1}\omega_2)=\alpha(\alpha(q_j,\omega_1),\omega_2)$$

For any word of length $n$, $\alpha(q_j,\omega_{1}\dots\omega_{n})$ is defined by inductively applying the lines above. We say a word is accepted if the state produced by the transition function is in $F$.

We will now define another type of finite state automata that will allow us to take an input alphabet $\Lambda$ to an output alphabet $\Gamma$.

\begin{definition}\label{ref:Mealy}
A Mealy machine is a tuple $M=(\Omega,\Lambda,\Gamma,\alpha,\beta,\omega_0,F)$, with a finite set $\Gamma$ known as the output alphabet and a function $\beta:\Omega\times \Lambda\to \Gamma$, which takes an input state and character as arguments and returns an output character. The other elements of the tuple are the same components of a finite state automaton.
\end{definition}

Similarly to $\alpha$, $\beta$ is defined inductively.

Due to the nature of the geometric addresses, it might not be clear how one would translate a particular geometric address into a dynamical address. Suppose that one selects a point that lives at the top of the basilica. The only partition this point could move to under the dynamics of the basilica is the left partition. However, the left partition maps to a union of the remaining three under the dynamics of $f$. The geometric address alone does not provide sufficient information on how a point's geometric address should translate into its dynamical one with the current collection of symbols. However, one can circumvent this issue by allowing the machine to output a word of one or two symbols at each step: from the top or bottom partition, the forced symbol $L$ is written together with the symbol that follows it. Each state and input still determine exactly one output word and one next state, so this solution lets us avoid a nondeterministic machine.

There is one more problem to address before we may construct our automaton. The current definition above only permits words of finite length to be read by the automaton. To address this, we need only extend our Mealy machine to an $\omega$-automaton, which may be defined as follows.

\begin{definition}\label{ref: Muller}
    A deterministic $\omega$-automaton is a tuple $M_{\omega}=(Q,\Lambda,\alpha,q_0,A)$, where the first four components are the same as a deterministic finite-state automaton and $A$ is the acceptance component. We say such an automaton is a Muller-automaton if $A\subset \mathscr{P}(Q)$.
\end{definition}

\begin{remark}
    Strictly speaking, the transition function is unchanged: $\alpha$ is still a map $Q\times\Lambda\to Q$, and only the acceptance component $A$ lives in the power set of $Q$. In particular, we may treat $\alpha$ as the FSA's transition function. 
\end{remark}

\begin{definition}\label{ref: Read}
    Let $M_{\omega}=(Q,\Lambda,\alpha,q_0,A)$ be a deterministic $\omega$-automaton. A run on a word $w\in \Lambda^{\omega}$ is a sequence $\rho\in Q^{\omega}$ such that

    $$\rho(j)=\alpha(\rho(j-1),w(j))$$

    for every $j\in\N$.
    
\end{definition}

Our next objective, by introducing the following definitions, is to extend the notion of acceptance to the case of infinite strings.

\begin{definition}\label{ref: Occurence and Infimum}
    For $w\in\Lambda^{\omega}$ (or, more generally, for any sequence over a finite set, such as a run $\rho\in Q^{\omega}$) we define the occurrence set and the infinity set as follows:

    $$\text{Occ(w)}=\{\lambda\in \Lambda: \exists j\in\N \ \text{s.t.} \ w(j)=\lambda\}$$
    
    $$\text{Inf(w)}=\{\lambda\in \Lambda: \forall i\exists j\in\N, j>i \ \text{s.t.} \ w(j)=\lambda\}$$
\end{definition}

\begin{remark}
   $\text{Inf}(w)$ is guaranteed to be finite as the alphabet is a finite set. This is because a word cannot consist of characters that exceed the alphabet; however, a word does not need to contain every character in the alphabet to be a member.
\end{remark}

With this definition in hand, we can define what acceptance means for our purposes:

\begin{definition}\label{ref: Acceptance of an infinite word}
    Let $M$ be an $\omega$-automaton. We say a word $w\in \Lambda^{\omega}$ is accepted if $\text{Inf}(\rho)\in A$, where $\rho$ is the run on $w$ and $A\subset\mathscr{P}(Q)$ is the acceptance component.
\end{definition}

With the following definitions in mind, we may now formally construct our Mealy machine $M_D$. This machine will consist of five different states: an initial state $i$ together with the four states $t,l,b,r$. The ambiguity that occurs when one has a point in the top partition and the bottom partition of the basilica is handled by output words of length two.

\newpage

Let us construct $M_{D}=(Q,\Lambda,\Gamma,q_0,\alpha,\beta,F)$ as follows:

$$Q=\{i,t,l,b,r\}$$
$$F=Q\setminus\{i\}$$
$$q_0=i$$
$$\Lambda=\{t,l,b,r,1,2,3\}$$
$$\Gamma=\{T,L,B,R\}$$
$$\alpha:Q\times \Lambda\to Q$$
$$\beta:Q\times \Lambda\to \Gamma^{+}$$
Here $\Gamma^{+}$ denotes the set of nonempty finite words over $\Gamma$; this is the only change from Definition~\ref{ref:Mealy}. The functions $\alpha$ and $\beta$ will be required to match the tables within the given input but are allowed to behave differently for any other given tuples, as these tuples will not be reachable by words in $\Lambda^{*}$. $M_{D}$ needs only to correctly translate words produced by the input alphabet to the output alphabet.

\begin{table}[h]
    \centering
    \begin{minipage}[t]{0.48\textwidth}
        \vspace{0pt}
        \centering
        \begin{tabular}{|c|c|c|}
            \hline
            State & Character & Output \\ \hline
            $i$ & $t$ & $t$ \\ \hline
            $i$ & $l$ & $l$ \\ \hline
            $i$ & $b$ & $b$ \\ \hline
            $i$ & $r$ & $r$ \\ \hline
            $t$ & $1$ & $b$ \\ \hline
            $t$ & $2$ & $r$ \\ \hline
            $t$ & $3$ & $t$ \\ \hline
            $l$ & $1$ & $b$ \\ \hline
            $l$ & $2$ & $r$ \\ \hline
            $l$ & $3$ & $t$ \\ \hline
            $b$ & $1$ & $b$ \\ \hline
            $b$ & $2$ & $r$ \\ \hline
            $b$ & $3$ & $t$ \\ \hline
            $r$ & $1$ & $b$ \\ \hline
            $r$ & $2$ & $r$ \\ \hline
            $r$ & $3$ & $t$ \\ \hline
        \end{tabular}
        \caption*{$\alpha$}
    \end{minipage}
    \hfill
    \begin{minipage}[t]{0.48\textwidth}
        \vspace{0pt}
        \centering
        \begin{tabular}{|c|c|c|}
            \hline
            State & Character & Output \\ \hline
            $i$ & $t$ & $T$ \\ \hline
            $i$ & $l$ & $L$ \\ \hline
            $i$ & $b$ & $B$ \\ \hline
            $i$ & $r$ & $R$ \\ \hline
            $t$ & $1$ & $LB$ \\ \hline
            $t$ & $2$ & $LR$ \\ \hline
            $t$ & $3$ & $LT$ \\ \hline
            $l$ & $1$ & $B$ \\ \hline
            $l$ & $2$ & $R$ \\ \hline
            $l$ & $3$ & $T$ \\ \hline
            $b$ & $1$ & $LB$ \\ \hline
            $b$ & $2$ & $LR$ \\ \hline
            $b$ & $3$ & $LT$ \\ \hline
            $r$ & $1$ & $B$ \\ \hline
            $r$ & $2$ & $R$ \\ \hline
            $r$ & $3$ & $T$ \\ \hline
        \end{tabular}
        \caption*{$\beta$}
    \end{minipage}
    \caption{Transition tables for $\alpha$ and $\beta$.}
    \label{tab:alpha-beta}
\end{table}


\newpage

We will define the following two self-maps within the geometric and dynamical descriptions of the Basilica.

$$\psi_{f}: \Gamma_{D} \to \Gamma_{D}$$
$$\psi_{f}(a_1a_2a_3a_4\dots)=a_2a_3a_4\dots$$
Here $\Gamma_{G}=\{t,l,b,r\}\times\{1,2,3\}^{\N}$, with the product topology, denotes the space of geometric addresses (Definition~\ref{def:geometric-address}), and $\omega$ denotes an arbitrary tail in $\{1,2,3\}^{\N}$.
$$\varphi_{f}: \Gamma_{G} \to \Gamma_{G}$$


\begin{gather*}
     t\omega \xrightarrow{\varphi_{f}} l\omega, 
    r1\omega\xrightarrow{\varphi_{f}} b\omega, \\
    r2\omega\xrightarrow{\varphi_{f}} r\omega, 
    r3\omega\xrightarrow{\varphi_{f}} t\omega\\
    b\omega \xrightarrow{\varphi_{f}} l\omega, \ l1\omega\xrightarrow{\varphi_{f}} b\omega, \\
    l2\omega\xrightarrow{\varphi_{f}} r\omega, \
    l3\omega\xrightarrow{\varphi_{f}} t\omega  
\end{gather*}

\begin{center}
{\Large
\begin{tikzcd}[sep=huge]
J_f \arrow[r, "\text{f}"]                         & J_f                             \\
\Gamma_{D} \arrow[r, "\psi_f"] \arrow[u, "\pi_D"]        & \Gamma_{D} \arrow[u, "\pi_D"']  \\
\Gamma_{G} \arrow[r, "\varphi_{f}"'] \arrow[u, "M_{D}"] & \Gamma_{G} \arrow[u, "M_{D}"']
\end{tikzcd}
}

\end{center}

\begin{theorem}
The diagram above commutes.
\end{theorem}

\begin{remark}
    Throughout the proof, we will abuse notation and write the resulting output of the Mealy machine as $M_D(\omega)$.
\end{remark}

\begin{proof}

The first square commutes via the construction of $\pi_{D}$: for $\omega\in\Sigma_D$, $\tau$ maps the arcs of $\omega$ onto those of $\psi_f(\omega)$, so $\theta(\psi_f(\omega))=\tau(\theta(\omega))$ and $\pi_D(\psi_f(\omega))=\ell(2\theta(\omega))=f(\pi_D(\omega))$. Every output of $M_D$ is admissible by Theorem~\ref{thm:transition}, so $M_D$ takes values in $\Sigma_D$, where $\pi_D$ is defined. Thus, most of the work of this proof lies in showing that the proposed automata for geometric addresses and dynamical addresses make the bottom square commute. The number of cases that need to be checked is finite; any geometric address can only start with the following:

$$\{t1,t2,t3,l1,l2,l3,r1,r2,r3,b1,b2,b3\}$$

 However, we also know that $r$ and $\ell$ have (relatively) symmetric behavior and $b$ and $t$ have symmetric behavior under the transition function. This reduces the number of cases to be checked significantly. This implies that, for a sufficient proof, one must check fewer cases. We will first give an example of one such case, and then treat all cases at once.


\emph{Example case.}
Suppose $z_0\in J_f$ and $[z_0]_{G}:\N\to \Lambda^{*}$ (the corresponding geometric address) is defined as $[z_0]_{G}=t\omega_1\omega_2...$. The self-map $\varphi_f$ takes $[z_0]_{G}$ to the word $[z_0']_{G}=l\omega_1\omega_2...$. Suppose $\omega_1=1,\omega_2=2$. The automaton $M_D$ takes these geometric addresses to the following dynamical addresses.

$$M_{D}([z_0']_{G})=LBLR\beta(r,\omega_3)...=[z_0']_{D}$$

$$M_{D}([z_0]_{G})=TLBLR\beta(r,\omega_3)...=[z_0]_{D}$$

However, applying the shift map in the dynamical address space to the image of $[z_0]_{D}$ yields $\psi_f([z_0]_{D})=LBLR\beta(r,\omega_3)...=[z_0']_{D}$. This implies that

$$\psi_f\circ M_{D}([z_0]_{G})=M_{D}\circ\varphi_f([z_0]_{G}).$$

 


In general, write $X_1=B$, $X_2=R$, $X_3=T$ and $s_1=b$, $s_2=r$, $s_3=t$, and let $M_x(w)$ denote the output of $M_D$ started in state $x\in\{t,l,b,r\}$ on a tail $w\in\{1,2,3\}^{\N}$. Then $M_D(xw)=X\,M_x(w)$, where $X$ is the capital letter of $x$, and by Table~\ref{tab:alpha-beta}, for $k\in\{1,2,3\}$,
$$M_t(kw)=M_b(kw)=L\,X_k\,M_{s_k}(w),\qquad M_l(kw)=M_r(kw)=X_k\,M_{s_k}(w).$$
Hence
$$\psi_f M_D(tw)=M_t(w)=L\,M_l(w)=M_D(lw)=M_D(\varphi_f(tw)),$$
and likewise for $b$. For $x\in\{l,r\}$,
$$\psi_f M_D(xkw)=M_x(kw)=X_k\,M_{s_k}(w)=M_D(s_kw)=M_D(\varphi_f(xkw)).$$
These are all the cases in the definition of $\varphi_f$, so the bottom square commutes.

\end{proof}

\section{Conclusion}

By implementing methods from topology, complex analysis, and automata theory, we were able to relate the dynamical orbit of a point to its location on the basilica. In addition, the fact that the diagram commutes allows for a more combinatorial exploration of the basilica via the use of the geometric address of points. The initial success with the basilica might suggest that other families of connected Julia sets could be analyzed similarly. In particular, if a connected Julia set may be described in terms of finitely many connected components, then a combinatorial description via the geometric location of points might be viable.

\section{Acknowledgement and AI Use}
The authors thank Professor Shuyi Weng for guidance throughout this project. His guidance provided many insights and directions to explore for this project.\\\\
AI is only used for formatting, finding typos and reference searching in this paper.
\bibliographystyle{plainnat}

\end{document}